**What is new on computed Lorenz strange attractors: chaos or numerical errors?**

Lun-Shin Yao
Department of Mechanical and Aerospace Engineering
Arizona State University
Tempe, Arizona 85287

[2003.yao@asu.edu](mailto:2003.yao@asu.edu)

Abstract

Discrete numerical methods with finite time-steps represent a practical technique to solve initial-value problems involving nonlinear differential equations. These methods seem particularly useful to the study of chaos since no analytical chaotic solution is currently available. Using the well-known Lorenz equations as an example, it is demonstrated that numerically computed results and their associated statistical properties are time-step dependent. There are two reasons for this behavior. First, chaotic differential equations are unstable so that any small error is amplified exponentially near an unstable manifold. The more serious and lesser-known reason is that stable and unstable manifolds of singular points associated with differential equations can form virtual separatrices. The existence of a virtual separatrix presents the possibility of a computed trajectory actually "jumping" through it due to the finite time-steps of discrete numerical methods. Such behavior violates the uniqueness theory of differential equations and amplifies the numerical errors *explosively*. These reasons imply that, even if computed results are bounded, their independence on time-step should be established before accepting them as useful numerical approximations to the true solution of the differential equations. However, due to these exponential and explosive amplifications of numerical errors, no computed chaotic solutions of differential equations independent of integration-time step have been found. *Thus, reports of computed non-periodic solutions of chaotic differential equations are simply consequences of unstably amplified truncation errors, and are not approximate solutions of the associated differential equations.*

## 1. Introduction

In spite of numerous attempts, a convincing proof of the existence of the geometric Lorenz attractor for the Lorenz's differential equations was not achieved until recently. Tucker (2002) provided a solution to this problem, which is the 14th of the 18 challenging mathematical problems defined by Smale (1998). Viana (2000), in a review of the historical advancement of structural stability theories, explained the difficulties that Tucker had to overcome to complete his proof. The next step in the progression of this important research topic is the construction of an existing solution. Without proof of the existence of chaotic solutions of nonlinear differential equations, Tucker's solution is only valid for *algebraic mappings*. Equivalently, a non-periodic structure associated with Smale's horseshoe is for algebraic mappings only, as is explained in section 5 when we discuss the Lorenz attractors.

A possible approach to this issue is direct numerical integration, which is a powerful, popular, and convenient tool for solving initial-value problems for nonlinear differential equations arising in science and engineering. The use of such tools introduces truncation and rounding errors that often have a major impact on the quality of "computed results," which are the product of an application of a certain algorithm in a certain computing environment. *The major goal of this paper is to demonstrate why approximate chaotic solutions cannot be constructed by the current generation of such methods.* The important consequence is that the existence of chaotic solutions of differential equations has never been adequately demonstrated in spite of reports to the contrary.

In all specific cases reported in this paper, "computed results" fail to provide an approximate solution in any sense to the original initial-value problem. These ideas are studied for a particular chaotic system based on the well-known Lorenz equations for variables $x(t)$, $y(t)$, and $z(t)$. It is demonstrated that the numerical integration of this chaotic system is extremely sensitive to the integration time-step. Coupled with this behavior is the fact that different integration time-steps yield both "computed results" and computed statistical properties of the associated attractors that are dramatically different. These statistical properties can often be associated with physical quantities of interest in applications from science and engineering. *Of course, any "computed results" whose statistical properties are sensitive to integration time-steps are not useful.*

The Rössler equations (Stuart & Humphries, 1996, pp. 532-533) display a sensitivity of "computed results" to integration time-steps. Rössler himself noted, in the fall of 1975 (Rössler 2000, p. 213), the important role-played by integration time-steps in his system. However, there was no attempt to connect this sensitivity to the statistical properties of the "similar" attractors shown in Figure 7.1 of Stuart & Humphries (1996), and no explanations for these observations were provided. The sensitivity of "computed results" to integration time-steps has also been studied for the one-dimensional Kuramoto-Sivashinsky equation, a nonlinear partial differential equation (Yao 2007).

The sensitivity to initial conditions has been studied for hyperbolic systems (Bowen 1975) and is an active topic for nearly hyperbolic systems (Dawson, Grebogi, Sauer, & Yorke 1994; Viana 2000; Palis 2000; Morales, Pacifico & Pujals 2002). In particular, the methods of shadowing for hyperbolic systems have shown that trajectories may be *locally* sensitive to initial conditions while being *globally* insensitive since true trajectories with adjusted initial conditions exist. Such trajectories are called shadowing trajectories, and lie very close to the long-time computed trajectories. However, systems of differential equations arising from physical applications are not hyperbolic systems. If the attractor is *transitive* (ergodic), all trajectories are inside the attractor so they are generally *believed* to be solutions of the underlying differential equations whatever the initial conditions. Little is known about systems, which are not transitive (non-ergodic). Do & Lai (2004) have provided a comprehensive review of previous work, and discussed the fundamental dynamical process indicating that a long-time shadowing of non-hyperbolic systems is not possible.

Two important facts about the computation of numerical chaotic solutions of differential equations, which are not commonly known, are:

1. No computed chaotic solution of the Lorenz system, which is independent of the integration time-step, exists. The same conclusion can be extended to other chaotic situations, including the direct numerical simulation of turbulence through the Navier-Stokes equations.

2.  A sensitivity-to-initial-condition is frequently viewed as a *necessary* and *sufficient* condition for the existence of chaos by most readers.  However, this property is also noted in the solutions of all nonlinear differential equations when the values of their governing parameters are larger than appropriate critical values (Ghosh Moulic & Yao 1996; Yao & Ghosh Moulic 1994, 1995a, 1995b; Yao 1999, 2007, 2009).  Consequently, it cannot be argued that this sensitivity is a sufficient condition for chaos.  It is worthwhile to note that the sensitivity to initial conditions associated with a set of nonlinear differential equations is a reflection of a characteristic of a physical system; on the other hand, integration time-step is an *artificial* computational quantity.  That a discrete numerical computation must not be time-step dependent in order to be considered as an approximate solution was first put forward by Von Neumann (Teixeira, Reynolds & Judd, 2007; Lorenz 2008; Yao & Hughes 2008a, 2008b).

In this paper, I discuss the first of these issues in detail, showing that it is a consequence of unavoidable numerical errors, and briefly consider the second.  In the next section, I show that such "computed results" as well as a corresponding long-time averaged statistical correlation display a sensitive dependence on integration time-steps.  A systematic decrease in the magnitude of the time steps does not lead to a convergent pattern; rather irregularly fluctuating results are noted due to instabilities.  Differences between different samples can be small in a suitable sense, but never demonstrate convergence, indicating that the bounded "computed results" are contaminated by numerical errors.  Consequently, there is no guarantee that any of such "solutions" is close to the correct one.

The third section contains a demonstration of the exponential amplification of a small difference between two trajectories, which involve different integration time-steps, occurring when they move in the direction of an unstable manifold (the x-axis for the specific cases treated here).  When two different, as just mentioned, trajectories move in the direction of a stable manifold, their difference becomes smaller. This difference depends on the initial error introduced by the different time-steps, and prevents the determination of an approximate computed solution by any discretized numerical method.

In the fourth section, I show that a significant amplification of numerical errors occurs when a trajectory, in violation of theoretical expectations, jumps through a two-dimensional virtual separatrix. In contrast to the behavior noted in section 3, this behavior is independent of the differences induced by the integration time-steps before amplification.

Section 5 shows that "computed results" can be used to construct a strange attractor, even though they are time-step dependent. Each integration time-step generates its own algebraic mapping. For such algebraic mappings, the computed trajectory seems to visit the edges of the attractor less frequently. This observation agrees with the finding by Tucker (2002) that the stretching rate near the edge of the attractor is smaller than $\sqrt{2}$, the required minimum value to ensure its transitivity (Viana 2000). Numerical computed results do not demonstrated the property of topological transitivity so they are not ergodic. Amazingly, it can be shown by following Tucker's computer assisted method that attractors generated with different time-steps and contaminated by truncation errors satisfy the properties of Smale's horseshoe.

Conclusions are discussed in the final section.

2. **Time-Step Sensitive Numerical Solutions**

I use the Lorenz equations,

$$\dot{x} = -sx + sy,$$
$$\dot{y} = rx - y - xz,$$
$$\dot{z} = -bz + xy,$$

(1)

with the widely used values of the parameters, s = 10, r = 28, b = 8/3, as the basis for showing that the "computed results" are time-step dependent and do not converge for large time.  The initial conditions are x = 1, y = −1, z = 10.  All results presented are generated by an explicit, second-order accurate Adams-Bashforth method.  The time history of x is plotted in Figure 1 for three different time steps, clearly showing the divergence of the "computed results."  Similar behavior was noted with Adams-Bashforth methods up to the fifth order; an implicit Crank-Nicholson method; second-order and forth-order Runge-Kutta methods; adaptive methods; and compact time-difference schemes.  In no case was a convergent solution obtained for $t \geq 20$.

The trajectory for one of the cases of Figure 1 ($\Delta t = 0.0001$) is shown in three projections in Figure 2 for the first 300,000 computational time steps.  This figure is useful as a geometric aid to identifying the location (and reasons) for the observed divergence. Initially, the trajectory moves smoothly around the two singular points (reverse spiral), and from one singular point to the other.  However, at a certain time identified by an arrow in the figure, the computed trajectory diverges.  This behavior signals the "break-down" of the computation; it will be discussed in greater detail in Section 4.

The apparent random behavior in Figure 1 is consistent with the common expectation that it is impossible to repeat the time histories of chaotic different equations due to their extreme sensitivity.  On the other hand, in order to be useful in scientific or engineering applications, the statistical properties of a chaotic "computed result," which can be of physical significance, must *not* be sensitive to the integration time-step.  As an example, Figure 3 displays the time-averaged $L_2$ norms corresponding to the "computed results" of Figure 1,

$$E(t) = \frac{1}{t} \int_0^t \left| x^2 \right| dt .$$

They also strongly depend on the integration time-step.  Computations for much longer times than the ones shown in the figure reveal that the various E(t) continue to be dependent on the integration time-step and do not converge.  For long time, E(t) can be interpreted as the moment of inertia of the numerical Lorenz attractor about x = 0 since the attractor is assumed topologically transitive.  Since E(t) is an important statistical and geometric property of the attractor, different E(t) implies different attractors.  The geometric properties of the "computed results" will be examined in the next section to show that the sensitivity to integration time-steps is the consequence of repeated amplification of numerical errors.

3. **Exponential Growth of Numerical Errors**

The following discussion is based on the normal form of the Lorenz system (Tucker 2002; Viana 2000). With s = 10, b = 8/3, and r = 28, these equations are

$$\dot{x} = 11.8x - 0.29(x+y)z,$$
$$\dot{y} = -22.8y + 0.29(x+y)z,$$
$$\dot{z} = -2.67z + (x+y)(2.2x - 1.3y)$$

(2)

Their three equilibrium points are $(0, 0, 0)$ and $(\pm 5.5929, \ \pm 2.8981, 26.8698)$. The coefficients of the linear terms in (2) are eigenvalues and represent the growth rates of x, y, and z, respectively. The non-linear terms lead to energy transfers among x, y and z (Yao 1999). Far away from the three equilibrium points, the net effect is attracting since the sum of the eigenvalues is negative. Once a trajectory moves close to the three equilibrium points, it is trapped in a complex attractor due to the competition between attraction and repulsion of the three equilibrium points. This shows that the Lorenz attractor is inside a large attracting open set; hence, it is *robust* (Viana 2000).

The Euler method is used to integrate equations (2). The initial conditions are the same as those used in Section 2. An explicit, second-order-accurate Adams-Bashforth method was also used, and showed that the results were not dependent on these two "low-accuracy" numerical methods. These results, when properly interpreted, identify two mechanisms that contribute to the sensitivity of the "computed results" to the integration time-step. The first of these mechanisms will be discussed in the following material, the second in Section 4. This conclusion is insensitive to any particular numerical method as long as it involves truncation errors since the computations have been repeated with high-order and higher-level finite-difference methods, interval methods, and series methods.

The first mechanism is associated with the movement of the trajectory toward the z-axis (Figure 4). Equations (2) shows that, as the value of (x + y) becomes small, the value of the non-linear terms in the z equation decreases to near zero. This causes z to move exponentially toward $z = 0$, but the time required to reach $z = 0$ is infinite. Simultaneously, the value of x increases exponentially so that the trajectory turns sharply toward the direction of increasing x since the x-axis is the unstable manifold. In order to show the sensitivity to the integration time-steps, two different time steps, but the same initial conditions, were used to integrate equations (2). The results show that the errors accumulated before the trajectory reaches the z-axis are amplified exponentially along the x direction and cause substantial numerical errors in the computation. It is clear that this amplification begins near the z-axis. The differences between the two trajectories for different integration time-steps decrease as they approach the z axis (along the direction of the stable manifold), but start to increase exponentially as they turn toward the direction of the unstable manifold. This is a typical example for a positive Lyapunov exponent, and also is a hint about the existence of Smale's horseshoe. As the "computed results" repeatedly pass through the region just described, the corresponding trajectories move further apart. The study of whether such trajectories are shadowable for nearly-hyperbolic systems is a current research topic (Palis 2000; Morales, Pacifico & Pujals 2002; Tucker 2002). Even though it might be shadowable, *its statistical properties, which have practical interest, cannot be determined* (Dawson, Grebogi, Sauer & Yorke 1994). Why should one take the effort to solve chaotic differential equations when no statistical properties of the computed results can be determined?

4. **Explosive Amplification of Numerical Errors**

The second mechanism occurs close to the z-axis where the trajectory can turn in two opposite directions depending on whether the trajectory arrives at positive or negative x (Figure 5) since the z-axis is the intersection of stable and unstable manifolds. This means that a small numerical error can be "explosively" amplified. The breakdown of the computed results presented in Figure 2 belongs to this class. The reason for this "unshadowable" amplification of numerical errors is explained below.

It will be demonstrated in Section 5 that the trajectory frequently visits the neighborhood of the z-axis ($0 \leq z \leq 15$), where the values of x and y are small. It is clear from equations (2) that, if the trajectory starts on the z-axis, it will stay on it forever so that the z-axis is an invariant set for the saddle at the origin. A trajectory in the *inset* of a limit point will approach the limit point asymptotically. The inset of an attractor is called its basin. The separatrix is defined as the complement of the basins of attraction. The initial state of a trajectory must belong to a separatrix if its future ($\omega$) limit set is not an attractor. Therefore, a separatrix consists of the insets of the non-attractive (or exceptional) limit sets. An *actual* separatrix separates basins. However, if it does not actually separate basins, it is called a *virtual* separatrix (Abraham & Shaw 1992). A computed trajectory cannot penetrate a separatrix since that would violate the uniqueness theorem. Thus, a computed trajectory that jumps through a separatrix means that the "computed results" violate the differential equations.

The following linearized analysis shows that the inset of the saddle point at the origin near the z-axis is a two-dimensional surface that includes the z-axis. This inset is a virtual separatrix embedded in the attractor. The first-order linearized version of equations (2) for small x and y, that is, near the z-axis, are

$$\begin{aligned}
\dot{x} &= 11.8x - 0.29(x+y)z, \\
\dot{y} &= -22.8y + 0.29(x+y)z, \\
\dot{z} &= -2.67z
\end{aligned} \tag{3}$$

The solution for z is simply

$$z = z_0 e^{-2.67t}, \tag{4}$$

where $0 \leq z_0 \leq 15$ is the initial z location inside the attractor. The trajectories that reach the z-axis are found by setting

$$x = h(y, z) = A(z)\, y + O(y^2), \tag{5}$$

where $h(0, z) = 0$ and A(z) is a function to be determined. Equation (5) defines the local stable and unstable manifolds near the z-axis.  Using (5), the first two of equations (3) become

$$\dot{x} = \frac{dh}{dy}\dot{y}$$
$$= \frac{dh}{dy}\left[-22.8y + 0.29(h+y)z\right] \tag{6}$$
$$= 11.8h - 0.29(h+y)z,$$

where h can be determined by solving (6) with (4).  Comparison of the result with (5) provides

$$A = B\left[1 \pm \left(1 - B^{-2}\right)^{1/2}\right],$$
$$B = \frac{59.66}{z} - 1, \tag{7}$$

where the minus sign is for trajectories approaching the z-axis (stable manifold), and the plus sign is for the trajectories leaving it (unstable manifold).  This shows that the local stable manifold is a two-dimensional surface forming, with the z-axis, the inset of the saddle point at the origin.  Moreover, this result shows that the local stable and unstable manifolds approach the y-axis and x-axis, respectively, as z decreases to zero.  A trajectory approaching the z-axis along the stable manifold can only move away along one branch of the unstable manifold without jumping through the virtual separatrix.  The computed trajectories in figures 3, 4, and 5 display this behavior.

It is clear that the stable and unstable manifolds act as virtual separatrices and roughly divide the x-y plane ( *Poincare'* map) into four quadrants locally near the z-axis.  All meaningful trajectories should only travel in the first and third quadrants, but a *computed* trajectory may mistakenly move into the second and fourth quadrants, two forbidden zones, after jumping through the stable manifolds because of numerical errors introduced by finite integration time-steps, as shown in Figures 5c and 6c.  Such numerical errors substantially alter the shape of the attractor; this matter will be further discussed in the next section.  Once a computed trajectory moves into a forbidden zone, it can return to its "proper" track only at the beginning of a period of "winding" away from one of the two fixed points above the origin and by forming the "wing of a butterfly."  Dawson, Grebogi, Sauer & Yorke (1994) have pointed out that a continuous shadowing trajectory cannot exist for such a trajectory, that is, it is *unshadowable*.

The shortcoming of a discrete numerical method is that it cannot exactly reach a surface of zero thickness.  It is obvious that one of the two computed trajectories, shown in Figures 5 and 6, has passed through the two-dimensional inset of the saddle point at the origin, thereby violating the uniqueness theorem.  In Figure 6a and 6b, four slightly different integration time steps were used.  The corresponding computed trajectories moved closer to the z-axis within a circle of radius $10^{-10}$.  It is interesting to note that the two computed trajectories (cases A and B), which did not jump through the virtual separatrix, agree with each other, as do the two computed trajectories (cases C and D) that jump through the virtual separatrix and violate the differential equations.  However, note that these two sets of computed trajectories are substantially different from each other.

A commonly cited computational example in chaos involves two solutions of slightly different initial conditions that remain "close" for some time interval and then diverge suddenly.  In fact, this behavior is often believed to be a characteristic of chaos.  More properly, this phenomenon is actually due to the explosive amplification of numerical errors, and violation of the differential equations noted above.

Before closing this section, the essence of explosive amplification of truncation errors, which may be the origin of *homoclinic explosions*, is summarized in the following theorem:

**Theorem:**  Numerical errors can cause a chaotic trajectory of the Lorenz differential equations to penetrate a separatrix.  Since a pseudotrajectory of a chaotic system of non-linear differential equations can move very close to a separatrix, however small numerical errors introduced by discrete numerical methods can cause the pseudotrajectory to penetrate the separatrix.  This behaviour violates the uniqueness theorem; thus, the trajectory cannot be considered a solution of the Lorenz differential equations and is therefore unshadowable.

Proof: The λ-lemma (Palis & de Melo 1982) guarantees that a chaotic trajectory can move closer to the local intersection of stable and unstable manifolds (the z-axis for the Lorenz system) than any pre-assigned value. Consequently, the trajectory will travel through the separatrix unless there is zero truncation error.

5. **Lorenz Attractor**

The lack of convergence in the results of Figures 1 and 3 is, at first glance, unexpected, but is real. Attempts to ignore this behavior frequently rely on the following three commonly believed erroneous arguments. However, they cannot withstand careful scrutiny as the remarks provided below show.

*Argument 1*: Since a necessary property of chaos is the presence of a positive Liapunov exponent, or a positive nonlinear exponential growth-rate, the truncation error introduced by various numerical methods can be amplified exponentially. Hence, erroneous solutions develop differently due to different truncation errors. This is equivalent to saying that the finite-difference equations, which approximate the differential equations, are unstable. Thus, since convergence requires stability and consistency, convergent computed results are not achievable. Such unstable cases are shadowable, that is, they remain sufficiently close to the true trajectory with slightly different initial conditions. However, as demonstrated in this paper, Argument 1 is not valid uniformly in the entire geometric space. The breakdown in the numerical solutions for chaos shown in Figures 1 and 3 is sudden, explosive, and unshadowable, but it is not only due to the exponential growth of numerical errors associated with an unstable manifold. *Even if computed results are shadowable, their statistical properties cannot be determined, implying a useless computation!*

*Argument 2*: It is well known that chaotic solutions of differential equations are sensitive to initial conditions. The different truncation errors associated with different integration time-steps, in effect, lead to a series of modified initial conditions for later times.

Consequently, computed chaotic solutions are integration time-step dependent, and cannot be considered to be an approximate, in any sense, solution of the differential equations.

On the other hand, as demonstrated before ( Ghosh Moulic & Yao 1996; Yao & Ghosh Moulic 1994, 1995a, 1995b; Yao 1999, 2007), stable long-time numerical solutions for the Navier-Stokes equations and the one-dimensional Kuramoto-Sivashinsky equation are sensitive to initial conditions, but are also convergent and independent of the integration time-steps. This shows that a solution sensitive to initial conditions is not necessarily sensitive to integration time-steps.

*Argument 3*: It is commonly believed that the existence of an attractor guarantees the long-time correctness of numerical computations of chaos, irrespective of the numerical errors that are inevitably present in any computation. Such a concept has never been proved, but it is customarily used to support the belief that numerical errors do not invalidate particular computed chaotic results among the community working on numerical solutions of dynamic systems.

A reason, which often leads researchers to believe that any incorrect computed trajectory is acceptable as long as it resides in an attractor, is the attractor's property of being *robust*. Unfortunately, the true mathematical definition of *a robust attractor* is less dramatic and simply means that an attractor is included in a large attracting open set as stated in the section 3; thus, the existence of attractors does not make *incorrect computations become correct*! Furthermore, this argument is incorrect because a computation contaminated by numerical error can escape an existing correct attractor and create another attractor, which is associated with the incorrect numerical results. This will be discussed below.

The locations where numerical errors are amplified can be better discussed within the framework of a particular example, the Lorenz attractor. It should be emphasized that I do not have a method to explicitly compute the true Lorenz attractor due to unavoidable numerical errors. I can only determine an erroneous attractor as others have. A computed Lorenz attractor for $\Delta t = 10^{-5}$, and the initial condition (1, -1, 10) is used for the following discussion. The computation is carried out for $10^8$ time-steps, and recorded every 1,000 time-steps. The attractor is constructed using 100,000 points, admittedly insufficient, but there are limitations due to the speed of the available computer.

Thin slices of the computed attractor normal to the z-axis are plotted at four different z locations in Figure 7. A short curve above the attractor shows that the computed trajectory rapidly enters the attractor from its initial location. This is because the attractor is *robust*. The bottom of the attractor looks like a "thin sheet," with the z-axis embedded in it, as shown in Figure 7a. For the purpose of demonstration, an expanded cartoon of the computed thin-attractor section of Figure 7a appears in Figure 9. It shows that the size of the attractor section contaminated by numerical errors is twice as large as the correct one, and its shape is also quite different. This suggests that the effect of numerical errors is by no means small.

Moving to z = 17.9, the attractor starts to "split" near its center and the z-axis is no longer embedded in it, as shown in Figure 7b. Above this value of z, the linearized analysis of Section 4 is not valid. For even larger values of z, the attractor splits into two parts due to the attraction from the two equilibrium points. The two large dots in Figure 7c locate the equilibrium points that no computed trajectory can reach; hence, there are two "holes" in the computed attractor near the two equilibrium points. In Figures 7 and 8, it is clear that the two-dimensional inset of the saddle at the origin connects the two-dimensional outsets of the other two fixed points.

Higher up, the cross-section of the computed attractor shrinks and finally disappears for z > 40. Thin slices of the computed attractor normal to the x-axis are plotted in Figure 8. For small x, the attractor splits into two symmetric parts, which look very much like a "butterfly," as is well known, when viewed from other angles. The attractor is very thin due to the strong contraction of the Lorenz equations (Viana 2000; Tucker 2002).

A single simulation, which is not an acceptable solution according to the present results, used to construct a numerical Lorenz attractor indicates that the orbit is *dense*, and the computed attractor seems *transitive* and *indecomposable*. The plots of Figures 7 and 8 seem to show that the attractor is finite and closed; hence, it is *compact* and *invariant* for a given time step, convincing evidence that it is an attractor satisfying the properties of Smale's horseshoe! The plots also show that the computed trajectory visits the edge of the attractor less frequently than its interior, a minor weakness. Since this trajectory is sensitive to the initial condition, it is commonly called a *strange* attractor (Viana 2000; Guckenheimer & Holmes 1983). However, recall that computed attractors are also sensitive to integration time-steps in that *different time-steps result in different computed attractors*, which all satisfy the mathematical properties of Smale's horseshoe, as demonstrated by Tucker (2002). The conclusion is that none of the attractors generated by an algebraic mapping is an acceptable solution of the differential equation even though they all satisfy the property of Smale's horseshoe.

Furthermore, studies of multiple solutions of the Navier-Stokes equations (Ghosh Moulic & Yao 1996; Yao & Ghosh Moulic 1994, 1995a, 1995b; Yao 1999, 2007, 2009) indicate that initial conditions can determine completely different long-time development of flow patterns and/or the frequencies and wave numbers of their fluctuations. Those computational results, *convergent and independent of integration time-steps*, can be obtained only for unstable flows not too far from their critical states, but are *sensitive to initial conditions*. This implies that large number of attracting open sets exist for unstable flows. The open sets can be disjoint or overlapping. The phenomena are certainly complex and different from the description of the simple *structure of Smale's horseshoe. The relevance of his horseshoe to differential equations is an open question.*

No convergent computational results can be found when the Reynolds numbers are much larger than the corresponding critical Reynolds number for unstable flows, a class of fluid flows that include turbulent flows. After we attempted to compute many of these flows, the reason for the lack of success became clear: it is impossible to construct a *stable* discretized numerical method, which is required by one of Von Neumann's criteria for convergence, without introducing sufficient numerical dumping (numerical viscosity). Since the viscous effect is small for large Reynolds number, introducing too much numerical viscosity contravenes the *consistency* requirement of Von Neumann's convergent criterion.

## 6. Conclusion

It has been demonstrated that attempts to compute numerical solutions of the Lorenz equations and their associated statistical properties are contaminated by errors due to the use of a discrete numerical method and finite computer arithmetic. Similar behavior has been discovered for the Rössler equations (Stuart & Humphries, 1996; Rössler 2000) and a particular one-dimensional partial differential equation, the Kuramoto-Sivashinsky equation (Yao 2007). Reasons for this behavior have been advanced. They suggest that nonlinear differential equations are not hyperbolic systems since they have discrete singular points. Each singular point has its own stable and unstable manifolds, which may form one or more virtual separatrices. Truncation errors of numerical computation are amplified along the unstable-manifold direction; hence, they violate the Von Neumann stability requirement necessary to ensure convergent solutions. The existence of a virtual separatrix allows a computed trajectory to "jump" through it. Such behavior violates the differential equations. *Even in the presence of bounded "computed results," their convergence should be examined before accepting them as useful numerical approximations to the solution of the differential equations.* There is no *rigorous* mathematical theory or any existing evidence that supports any other conclusion.

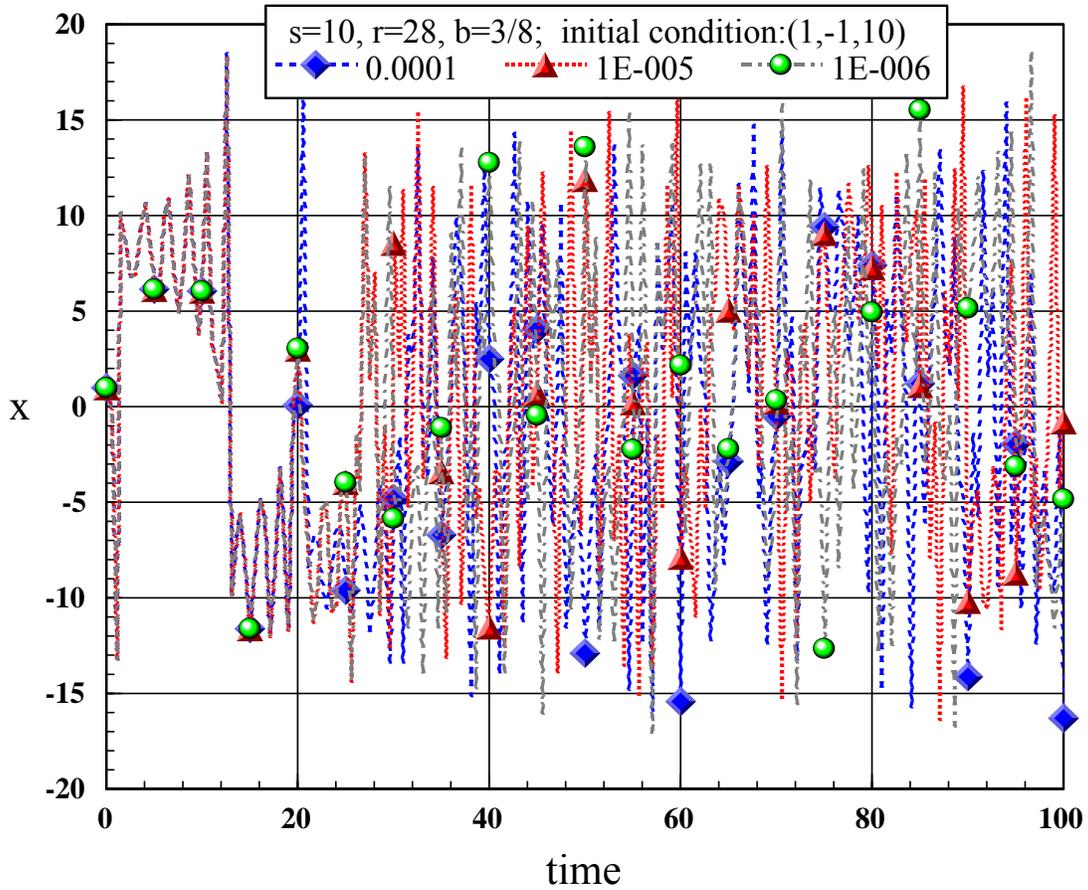

Figure 1.  Plots of the time history of x for three different time-steps
showing that the "computed results" for a time-step of 0.0001
start to diverge from the other two at approximately t=20.  They
all become different after t = 30.

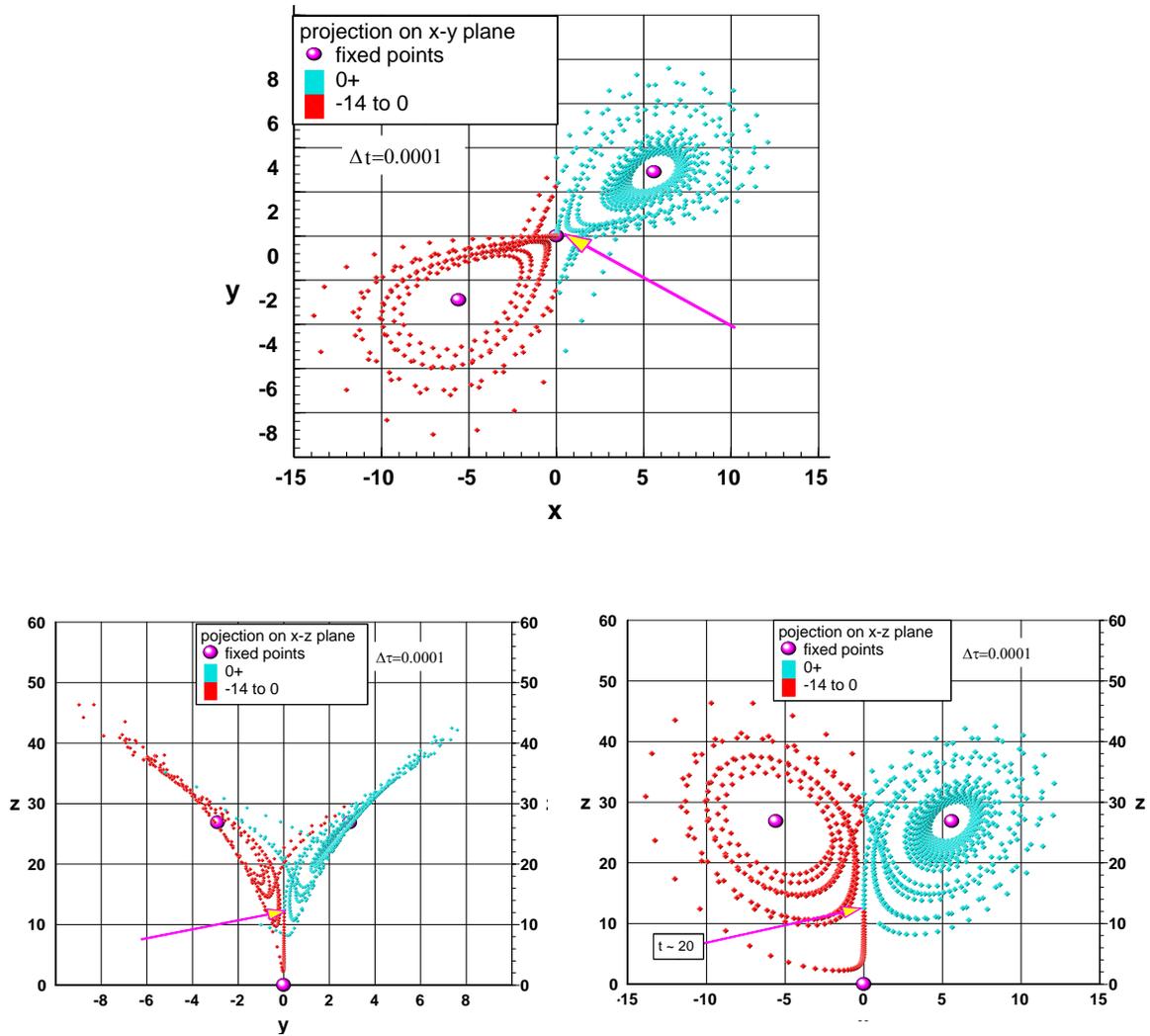

Figure 2. Plots of three projections for the first 300,000 computed points with a time-step of 0.0001. An arrow marks the location where the computed trajectory starts to deviate from those for other time-steps. This divergence occurs at time slightly smaller than 20.

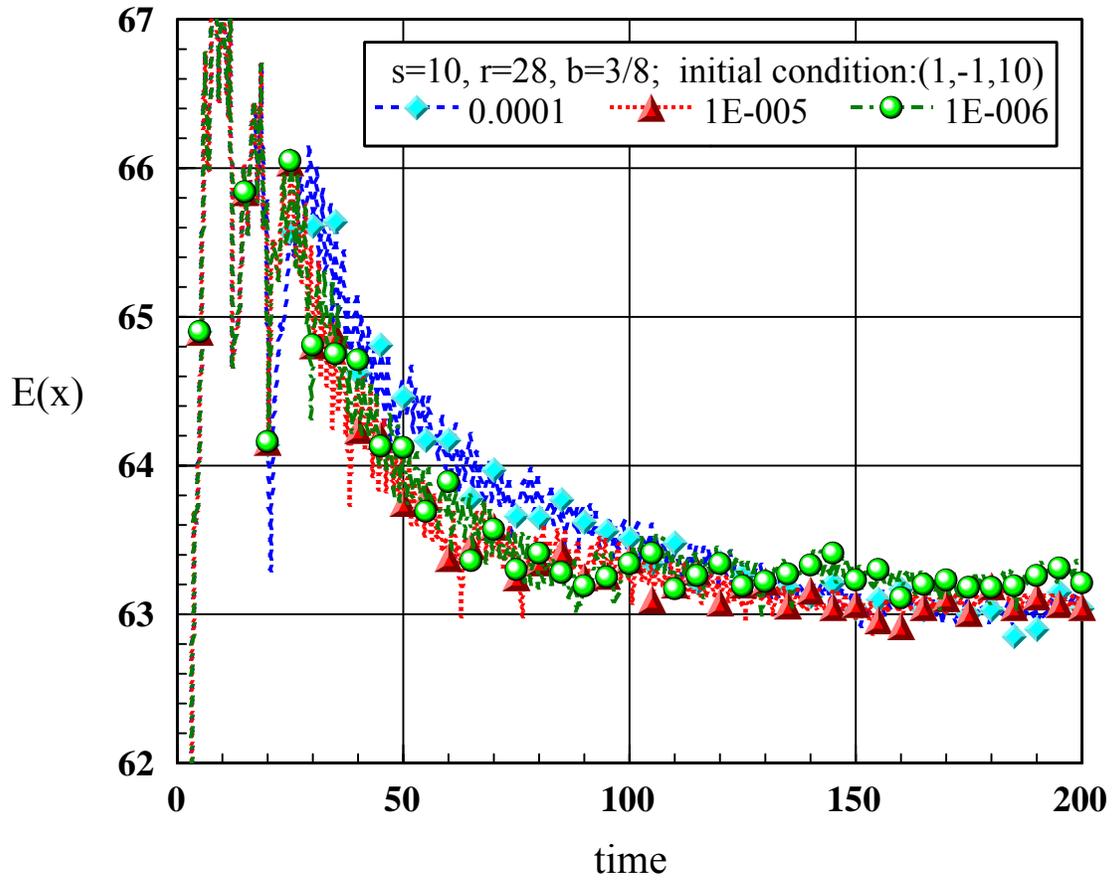

Figure 3.  The autocorrelations of computed x(t) (Section 2) for
three different time-steps showing that they are all
different.  Extending the computation to t=2000 did not
improve the convergence, and shows that E(x)
continuously fluctuates and does not show any
tendency to approach a constant.

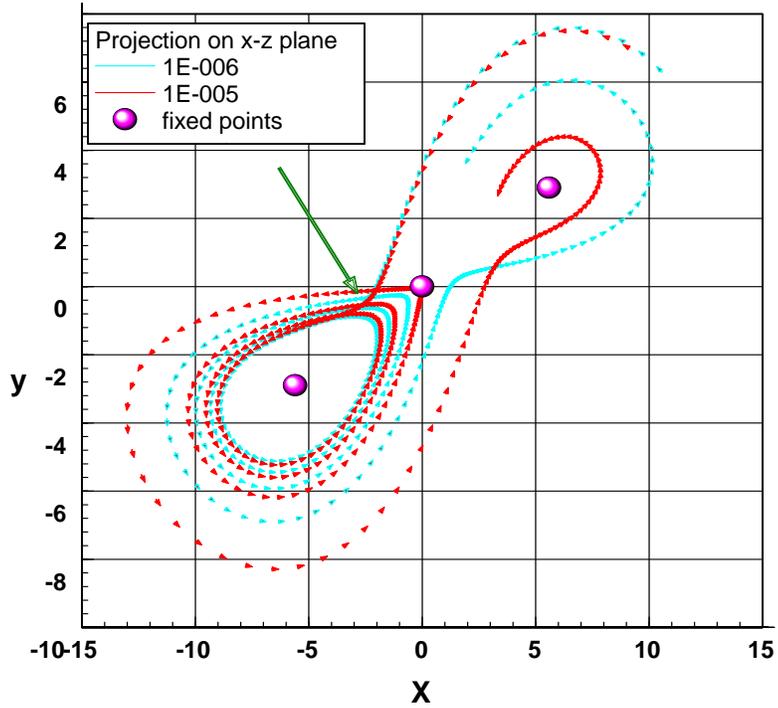

Figure 4a. "Computed results" for two different time-steps with identical initial conditions showing the exponential growth of the numerical errors. The arrow marks the starting location of dramatic error amplification.

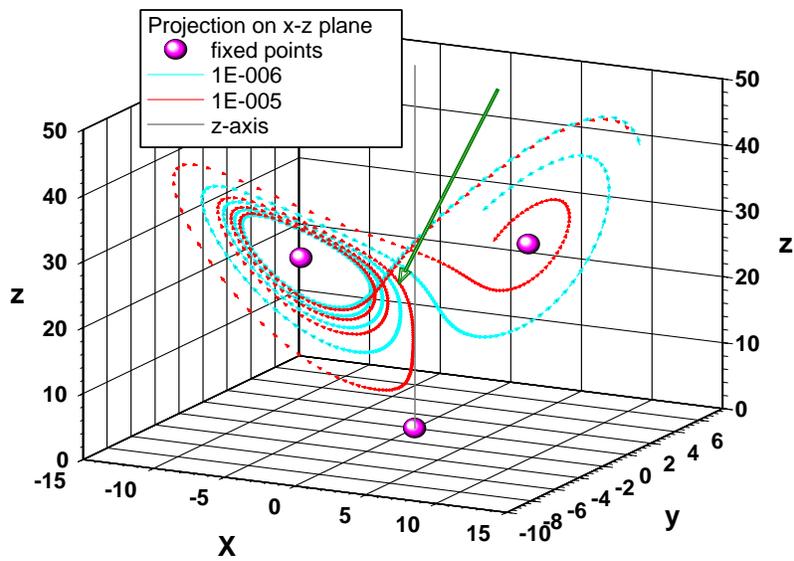

Figure 4b. A three-dimensional view of Figure 4a.

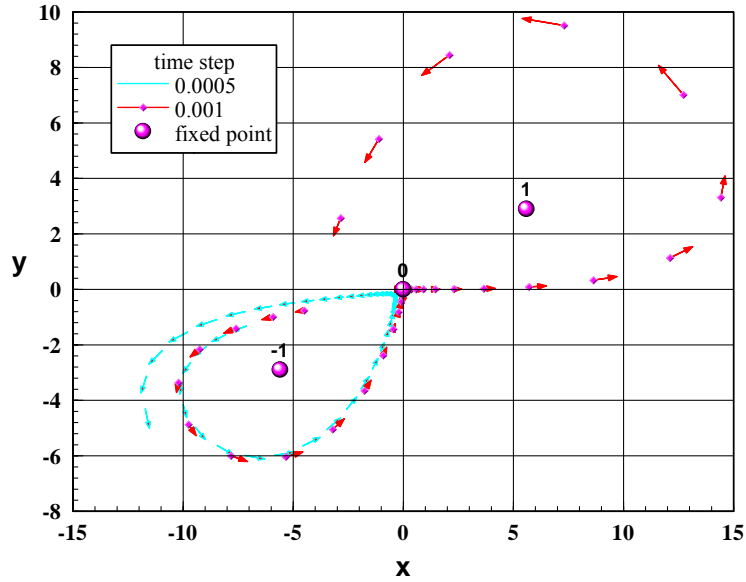

Figure 5a. "Computed results" for two different time-steps with identical initial conditions showing dramatic growth of the numerical errors. The trajectory associated with time-step of 0.001 jumps through the virtual separatrix.

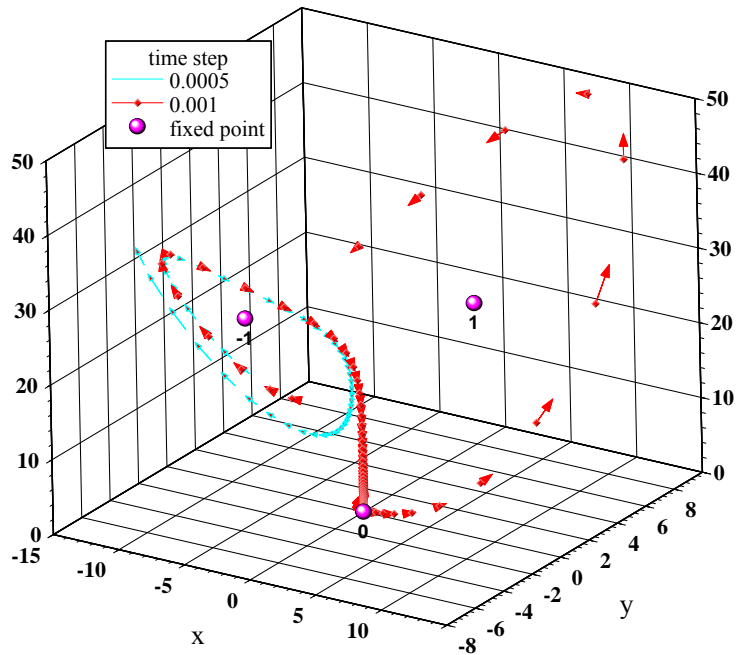

Figure 5b. A three-dimensional view of Figure 5a.

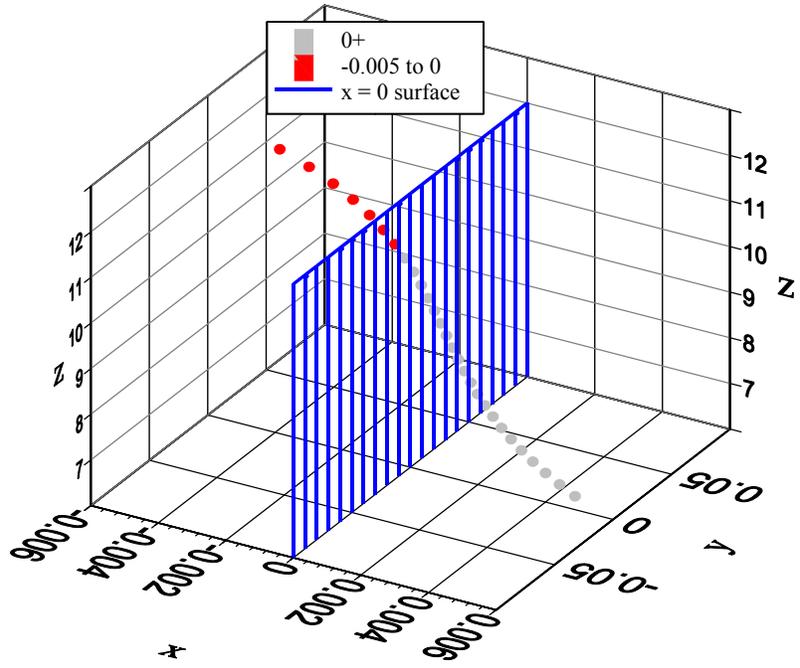

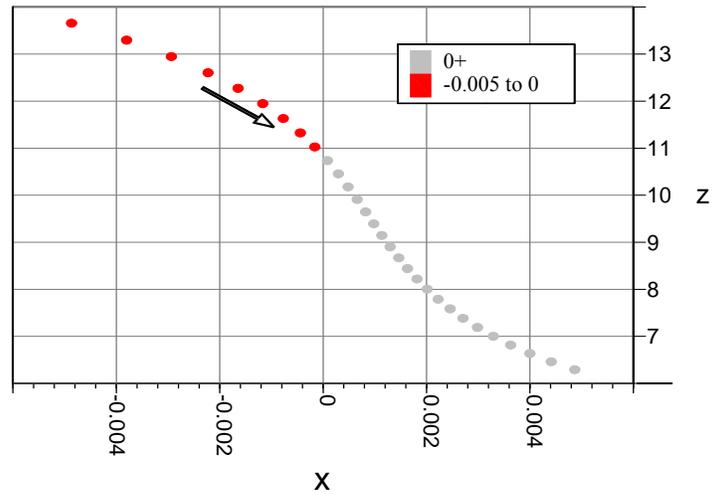

Figure 5c. Amplified views of Figure 5a & 5b.

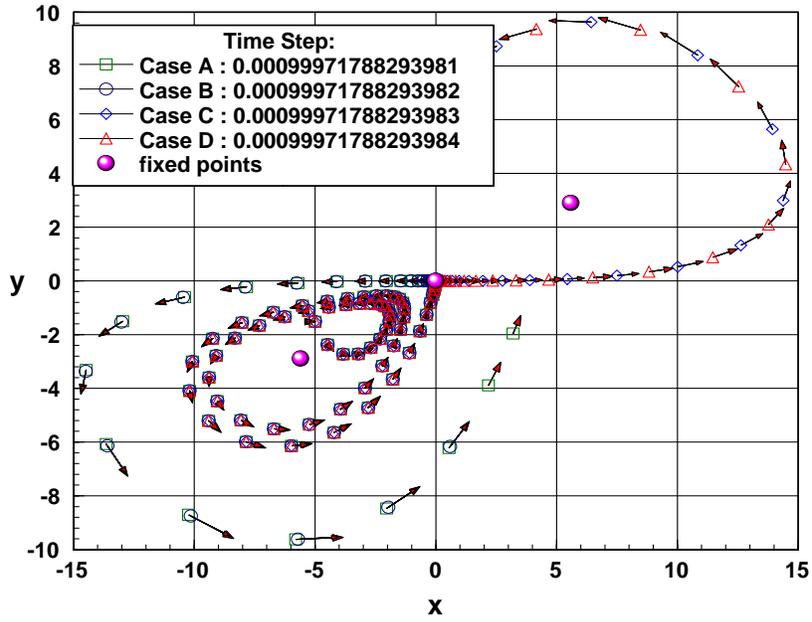

Figure 6a. "Computed results" showing that a dramatic growth of numerical errors occurs with time-steps of extremely small differences. Two of the computed trajectories jump through the virtual separatrix.

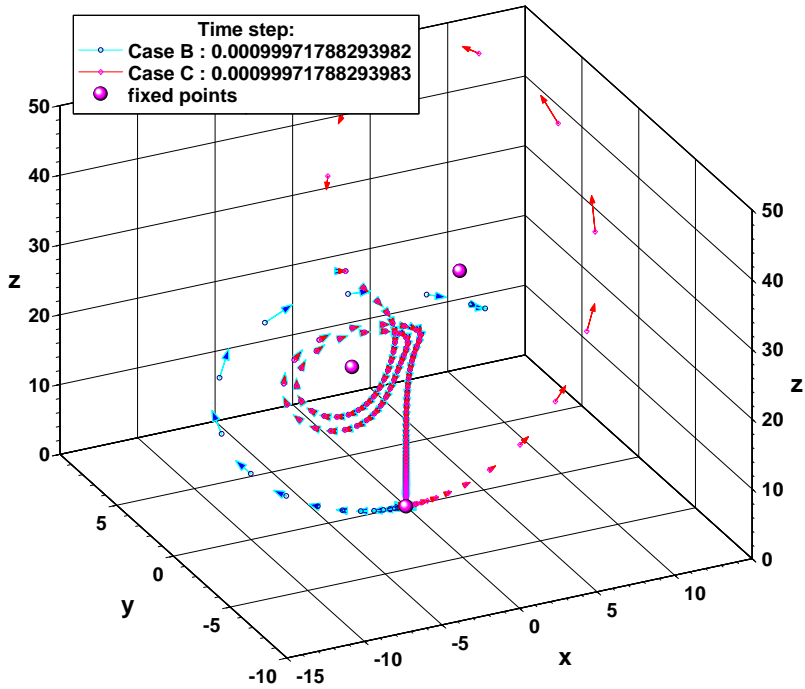

Figure 6b. A three-dimensional view of Figure 6a.

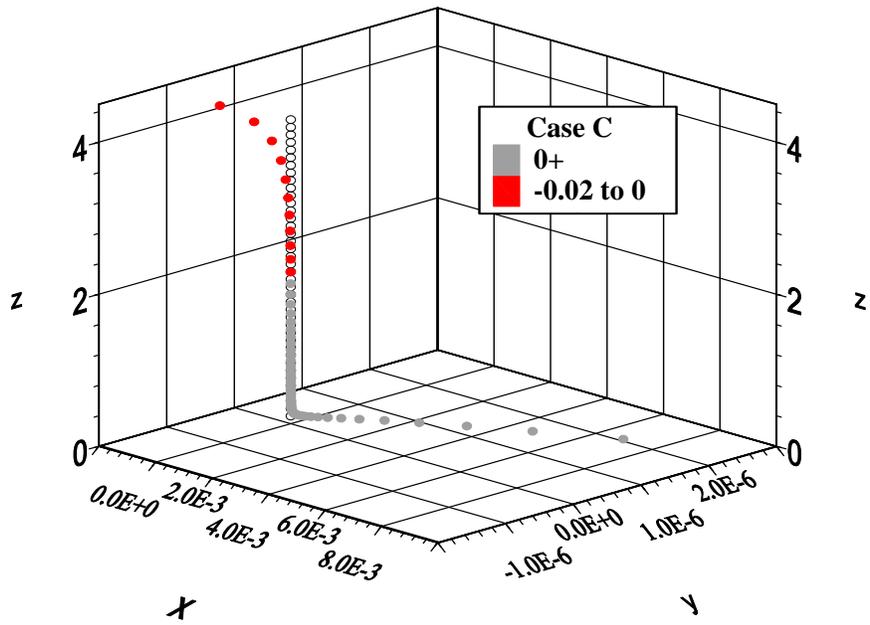

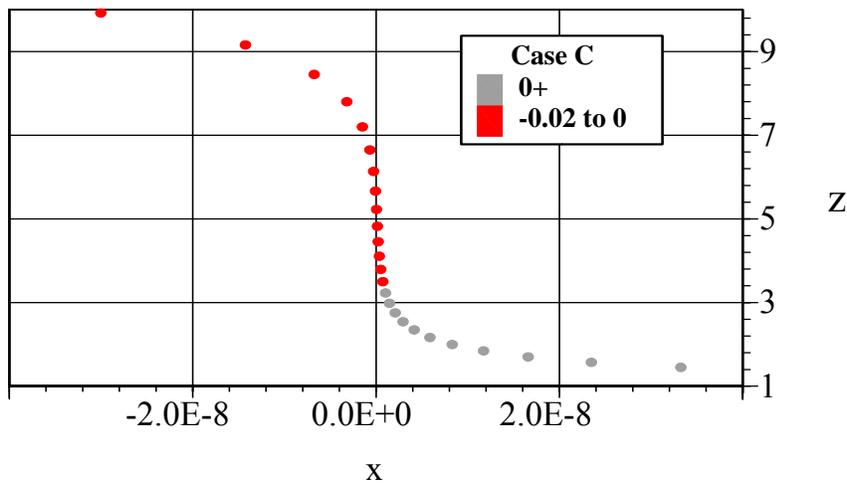

Figure 6c.  Amplified views of Figure 6a & 6b.

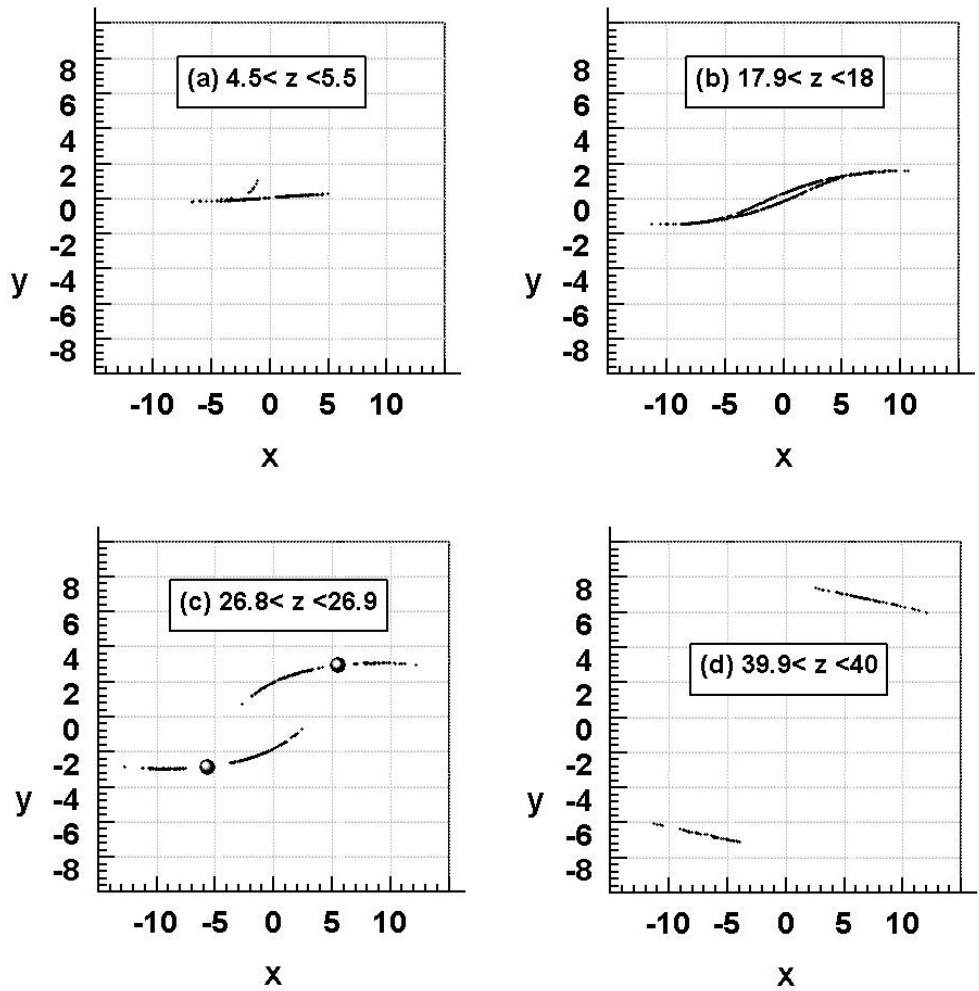

Figure 7.  Cross-sections of a computed Lorenz attractor at selected z locations showing that its thickness is thin and that its point density is not uniform.

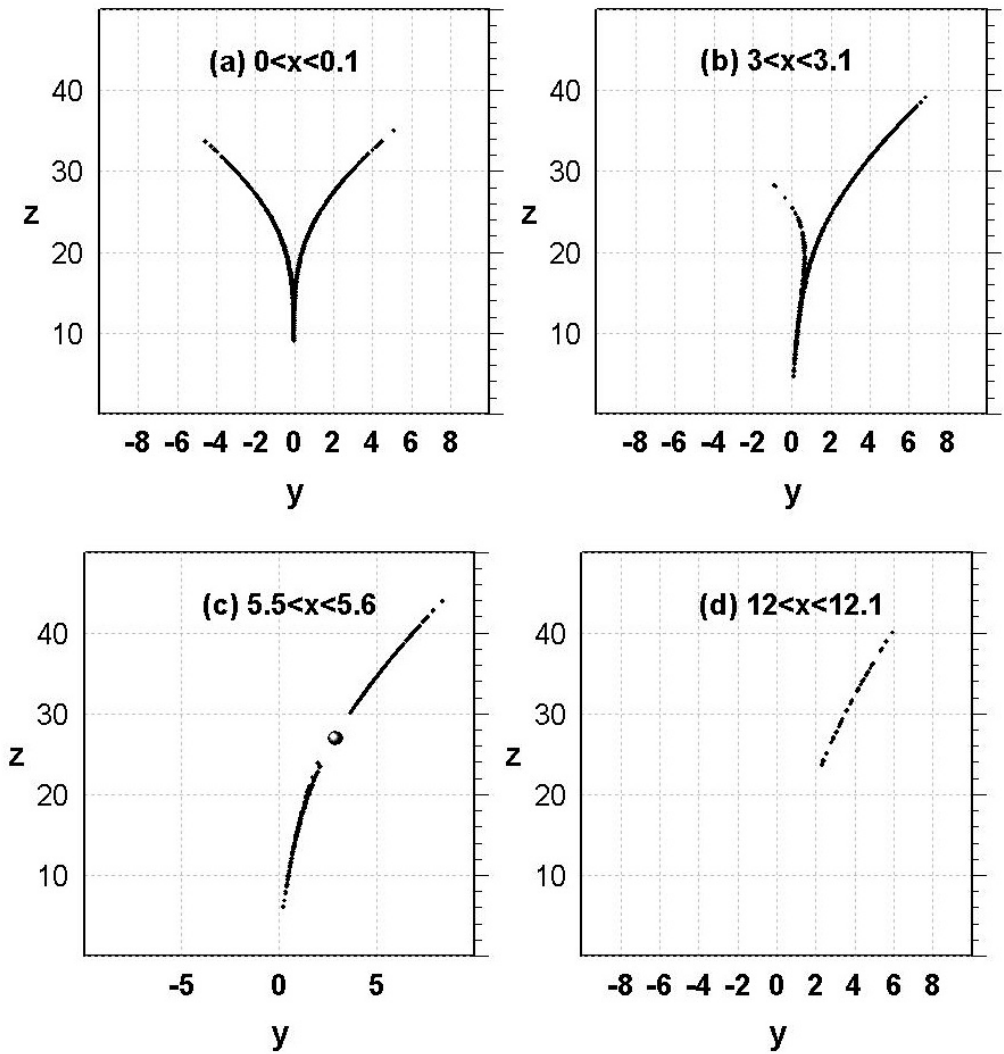

Figure 8. Cross-sections of a computed Lorenz attractor at selected x locations showing that its thickness is thin and that its point density is not uniform.

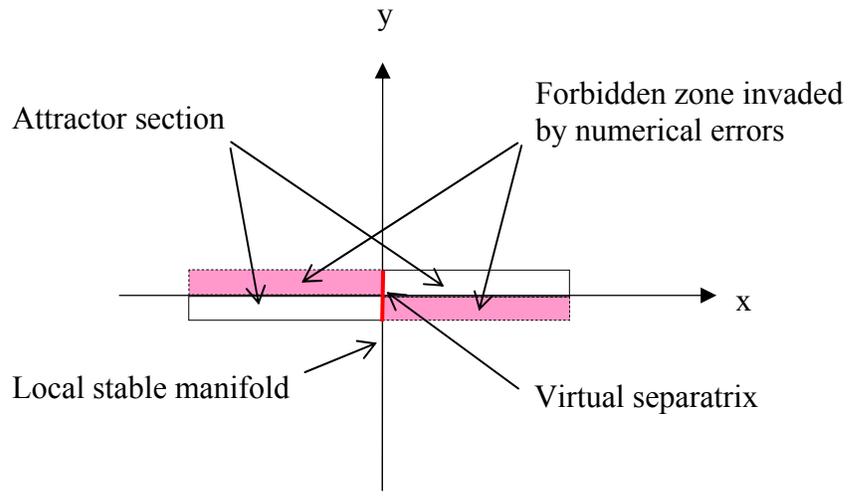

Figure 9.  Enlarged cartoon of Figure 7a after mapping the local unstable manifold to the x-axis.  The local stable manifold inside the attractor forms a finite virtual separatrix whose height is about $0 \le z \le 15$  (see equations 5 & 7).  The explosive error amplification causes the symmetry property of the computed results different from that of the Lorenz system, alters the shape of the Poincare section of the attractor, and doubles its area.